\documentclass{amsart}

\def\NN{\mathbb N}
\def\uu{\mathcal U}

\newcommand{\set}[1]{\left\lbrace #1\right\rbrace}

\newcommand{\qtq}[1]{\quad \text{#1}\quad }

\newtheorem{theorem}{Theorem}[section]
\newtheorem{proposition}[theorem]{Proposition}
\newtheorem{lemma}[theorem]{Lemma}
\newtheorem{corollary}[theorem]{Corollary}

\theoremstyle{definition}

\theoremstyle{remark}

\numberwithin{equation}{section}

\usepackage{amsmath,amsthm,amssymb}
\usepackage{enumerate}
\usepackage[]{lineno}
\usepackage{amsfonts}
\usepackage{latexsym}
\usepackage[T1]{fontenc}    
\usepackage[utf8]{inputenc}
\usepackage{enumerate}
\usepackage{color}
\begin{document}

\title[Multiple expansions]{Metric results for numbers with multiple $q$-expansions}
\author[S. Baker]{Simon Baker}
\address{School of Mathematics, University of Birmingham, Birmingham,  B15 2TT, UK}
\email{simonbaker412@gmail.com}
\author[Y. R. Zou]{Yuru Zou}
\address{College of Mathematics and Statistics,
Shenzhen University,
Shenzhen 518060,
People's Republic of China}
\email{yuruzou@szu.edu.cn}

\begin{abstract}
Let $M$ be a positive integer and $q\in (1, M+1]$.
A $q$-expansion of a real number $x$ is a sequence $(c_i)=c_1c_2\cdots$ with $c_i\in \{0,1,\ldots, M\}$ such that $x=\sum_{i=1}^{\infty}c_iq^{-i}$.
In this paper we study the set $\uu_q^j$ consisting of those real numbers having exactly $j$ $q$-expansions. Our main result is that
for Lebesgue almost every $q\in (q_{KL}, M+1),$ we have $$\dim_{H}\uu_{q}^{j}\leq \max\{0, 2\dim_H\uu_q-1\}\text{ for all } j\in\{2,3,\ldots\}.$$ Here $q_{KL}$ is the Komornik–Loreti constant. As a corollary of this result, we show that for any $j\in\{2,3,\ldots\},$ the function mapping $q$ to $\dim_{H}\uu_{q}^{j}$ is not continuous.

\end{abstract}

\date{\today}

\maketitle

\section{Introduction}
Fix a positive integer $M$. For $q\in (1, M+1]$ we call a sequence $(c_i)=c_1c_2\cdots\in \{0, 1,\cdots,M\}^{\mathbb{N}}$ a $q$-expansion of $x$ in base $q$ if
$$x=\pi_q((c_i)):=\sum_{i=1}^\infty \frac{c_i}{q_i}.$$
The study of $q$-expansions
 was pioneered in the papers of R\'enyi \cite{R1957} and Parry
\cite{Parry 1960}. Since these beginnings, the study of $q$-expansions has drawn significant attention. This is in part due to its connections with many other areas of mathematics. These areas include Dynamical Systems, Fractal Geometry, and Number Theory.

It is well-known that for $q\in(1, M+1],$  a number $x$ has an expansion in base $q$ if and only if $x\in I_q:=[0, M/(q-1)]$.
 When $q=M+1$ then  $I_q=[0,1]$ and every $x\in [0,1]$ has a unique expansion except for a countable set of exceptions that have precisely two.
When $q\in(1,M+1)$ then the situation is much more interesting. For instance, for any $q\in(1,M+1),$ it is the case that Lebesgue almost every $x\in I_q$ has a continuum of $q$-expansions \cite{S2003,DD2007}. For completion we mention that if $q>M+1$ then the set of points with an expansion in base $q$ is a fractal set, and if a point has an expansion in base $q$ then this sequence must be unique.

For each $j\in\mathbb{N}$ let
$$\uu_q^j:=\left\{x\in I_q: \#\pi_{q}^{-1}(x)=j \right\}.$$ Similarly, we let $$\uu_q^{\aleph_0}:=\left\{x\in I_q: \pi_{q}^{-1}(x) \text{ is an infinite countable set} \right\}.$$
We call $\uu_q^1$ the \emph{univoque set}. For simplicity, throughout this paper we write $\uu_q$ instead of $\uu_q^1$. We also let $U_q:=\pi_q^{-1}(\uu_q)$ be the corresponding set of sequences. The sets $\uu_{q}$ and $U_q$ were first properly studied by Erd\H{o}s et al. in the early 1990s \cite{EJK1990, EHJ1991, EJ1992}. Since then these sets have received significant attention and we now have a good understanding of their combinatorial, topological and fractal properties. See for instance the papers  \cite{DK1995,
Kallos1999,
Kallos2001,
GS2001,
DK2009,
KLD2010,
KL2015,
KKL2017,ABBK2019} and the results therein.
For $j\in\{2,3,\ldots\}\cup\{\aleph_0\}$, many important theorems have been obtained on the properties of the set $\uu_{q}^j,$ see
\cite{EJ1992,
S2009,
BS2014,
B2014,
B2015,
KLZ2017,
ZK2015,
ZWLB2016,
KK2018}.
However our knowledge of the set $\uu_{q}^j$ is significantly less than that of the set $\uu_{q}$. This paper is in part motivated by a desire to address this shortcoming.



Very little is known about the metric properties of $\uu_q^j.$ A simple bifurcation argument of Sidorov \cite{S2009}
and the first author \cite{B2015} implies that
\begin{equation}\label{e11qq}
\dim_H\uu_q^j\le{\dim_H\uu_q}\qtq{for all}q>1\qtq{and}
j\in\{2,3,\ldots\}\cup\{\aleph_0\}.
\end{equation} See also \cite[Lemma 5.5, Proposition 5.6]{ZKL2021}. Here and throughout $\dim_HF$ denotes the Hausdorff dimension of a set $F$. When $q$ is such that $\dim_H \uu_q=0$ then \eqref{e11qq} immediately implies that
\begin{equation}\label{e12}
\dim_H\uu_q^j={\dim_H\uu_q}\qtq{for all}
j\in\{2,3,\ldots\}\cup\{\aleph_0\}.
\end{equation}
Glendinning and Sidorov \cite{GS2001} showed  that  $\dim_H\uu_q=0$
for all $q\in(1, q_{KL}]$, where
$q_{KL}$ is the Komornik–Loreti constant (see Section \ref{s2} for more details). Therefore \eqref{e12} holds for all $q\in (1,q_{KL}].$ On the other hand, when $q=M+1$,
$\uu_q^2$ is countable and $\uu_q=[0,1]\setminus\uu_q^2$, so $\dim_H\uu_q^2=0<1=\dim_H\uu_q$.
When $q>M+1$, because every expansion is unique we have $1=\dim_{H}\uu_{q}>\dim_H\uu_q^j=0$. Therefore \eqref{e12} fails for every $q\geq M+1$. Because of these observations, it is natural to restrict our attention to studying the metric properties of $\uu_{q}^j$ for $q$ in the interval $(q_{KL}, M+1).$ Recently 
Sidorov \cite{S2009} showed that if $M=1$ and $q\approx 1.83929$ is the \emph{Tribonacci number}, i.e., the positive root of the equation $x^3=x^2+x+1$, then $\dim_H\uu_q^j={\dim_H\uu_q}$ for all $j\in \{2,3,\ldots\}$.
Motivated by Sidorov's work, the second author and her coauthors proved in \cite{ZKL2021} that there exist
infinitely many $q\in (q_{KL}, M+1)$
such that $\dim_H\uu_q^j={\dim_H\uu_q}$ for all $j\in \{2,3,\ldots\}$. In this paper we prove that for a typical $q\in (q_{KL},M+1)$ property \eqref{e12} does not hold. In particular we show that:
\begin{theorem}\label{t11}
For Lebesgue almost every $q\in (q_{KL}, M+1)$ we have $$\dim_{H}\uu_{q}^{j}\leq \max\{0, 2\dim_H\uu_q-1\}\text{ for all } j\in \{2,3,\ldots\}.$$  
\end{theorem}
For any $q\in (1,M+1)$ it is the case that $\dim_{H}\uu_{q}<1.$ Therefore Theorem \ref{t11} implies that $\dim_{H}\uu_{q}^{j}<\dim_{H}\uu_q$ for Lebesgue almost every $q\in(q_{KL},M+1).$ Notice that as well as establishing that $\dim_{H}\uu_{q}^{j}<\dim_{H}\uu_q$ holds almost everywhere in $(q_{KL},M+1)$, Theorem \ref{t11} also provides a lower bound for the difference between these quantities. The following theorem provides information on when $\uu_{q}^j$ is empty.
\begin{theorem}\label{t12}
Let $$ O:=\left\{q\in (q_{KL}, M+1):\dim_{H}\uu_q<1/2\right\}.$$ Then for Lebesgue almost every $q\in O$ we have $\uu_{q}^{j}=\emptyset$ for any $j\in \{2,3,\ldots\}$.
\end{theorem}The following corollary is a consequence of Theorem \ref{t11} and the results from \cite{ZKL2021} mentioned above. 

\begin{corollary}
	\label{continuity}
For any $j\in \{2,3,\ldots\},$ the function $f:(q_{KL},M+1)\to [0,1]$ given by $f(q)=\dim_{H}\uu_{q}^{j}$ is not continuous.
\end{corollary}Corollary \ref{continuity} is contrary to the case when $j=1$ for which it is known that the function mapping $q$ to $\dim_{H}\uu_{q}$ is continuous, see \cite{KKL2017,AK2019}. 

The rest of the paper is arranged as follows. In Section 2 we recall some relevant definitions and results from expansions in non-integer bases. In Section 3 we prove a number of technical results that will assist in our proof of Theorems \ref{t11} and \ref{t12}. In Section 4
we prove Theorems \ref{t11} and \ref{t12}.

\section{ Preliminaries}\label{s2}

Fix a positive integer $M$. We will denote an element of $\set{0,\ldots,M}^{\NN}$ by $(c_i)$ or $c_1c_2\cdots$.
We call a finite string of digits $w=c_1\cdots c_n$ with $c_i\in \set{0,1,\ldots, M}$ a \emph{word}. For convenience we let $\{0,\ldots,M\}^0$ denote the set consisting of the empty word.
Given two finite words
$w=c_1\cdots c_n$ and $v=d_1\cdots d_m$, we
denote by $wv=c_1\cdots c_n d_1\cdots d_m$ their concatenation.
Accordingly, for $k\in\mathbb{N}$ and a finite word $w$, we denote by $w^k$ or $(w)^k$
the concatenation of $w$ with itself $k$ times, and by $w^\infty$ or $(w)^\infty$ the concatenation of $w$ with itself infinitely many times. For a sequence $(c_i)$
 we denote by $\overline{(c_i)}=(M-c_1)(M-c_2)\cdots$
its reflection. 

We will use the lexicographic ordering on sequences.
 If $(c_i)$ and $(d_i)$ are two sequences, then we write $(c_i)\prec(d_i)$  if there exists $k\in\NN$ such that $c_i=d_i$ for $i=1,\ldots,k-1$ and $c_k<d_k$. Similarly we write $(c_i)\preceq (d_i)$  if $(c_i)\prec(d_i)$ or $(c_i)=(d_i)$.
We also write  $(d_i)\succ(c_i)$ if $(c_i)\prec(d_i)$, and $(d_i)\succeq (c_i)$ if $(c_i)\preceq (d_i)$.

For any fixed base $q\in(1,M+1]$, every $x\in I_q$ has a lexicographically largest expansion $b(x,q)=(b_i)$ obtained by the greedy algorithm, and a lexicographically largest infinite expansion $a(x,q)=(a_i)$; see \cite{BK2007,DKL2016}.
Such expansions are called the \emph{greedy} and \emph{quasi-greedy expansions} of $x$ in base $q$, respectively.
 A sequence $(c_i)$ is called \emph{finite} if it has a last non-zero digit, and it is called \emph{infinite} otherwise.
The case $x=1$ is particularly important. In this special case we introduce the simpler notation $\alpha(q):=a(1,q)=(\alpha_i)$.

We recall from \cite{KL1998, KL2002} that there exists a smallest base $q_{KL}\in (1,M+1)$ (depending
on $M$) for which $x=1$ has a unique expansion. This number is called the Komornik–Loreti constant and is defined using the classical Thue-Morse
sequence $(\tau_i)_{i=0}^\infty$:
\begin{equation}\label{e33qq}
\alpha(q_{KL}):=
\begin{cases}
(k+\tau_i)_{i=1}^{\infty} &\text{if $M=2k+1$,}\\
(k+\tau_i-\tau_{i-1})_{i=1}^{\infty}&\text{if $M=2k$,}
\end{cases}
\end{equation}
where $(\tau_i)_{i=1}^\infty=1101\ 0011\ \cdots$ denotes the truncated Thue--Morse sequence. Then the sequence $\alpha(q_{KL})$ begins with
\begin{equation}\label{12}
\begin{cases}
(k+1)(k+1)k(k+1)kk(k+1)(k+1)\cdots&\text{if $M=2k+1$,}\\
(k+1)k(k-1)(k+1)(k-1)k(k+1)k\cdots &\text{if $M=2k$.}
\end{cases}
\end{equation}

The following  lexicographic characterization of the  quasi-greedy  and greedy expansions was given in \cite{BK2007}.

\begin{lemma}
	\label{quasi-greedy lemma}
The following statements are true:
\begin{enumerate}[\upshape (i)]
 \item The map $q\mapsto \alpha(q)$ is a strictly increasing bijection between the interval $(1,M+1]$ and the set of infinite sequences $(\alpha_i)$ satisfying the lexicographic inequalities
\begin{equation*}
(\alpha_{n+i})\preceq (\alpha_i)\qtq{for all } n\geq 0.
\end{equation*}
\item For a fixed $q\in(1,M+1],$ the map $x\mapsto a(x,q)$ is a strictly increasing bijection between the interval $(0, M/(q-1)]$ and the infinite sequences $(a_i)$ satisfying the lexicographic inequalities
\begin{equation*}
(a_{n+i})\preceq \alpha(q)\qtq{whenever}a_n<M.
\end{equation*}
\item For a fixed $q\in(1,M+1],$ the map $x\mapsto b(x,q)$ is a strictly increasing bijection between the interval $[0, M/(q-1)]$ and the sequences $(b_i)$ satisfying the lexicographic inequalities
\begin{equation*}
(b_{n+i})\prec\alpha(q)\qtq{whenever} b_n<M.
\end{equation*}
\end{enumerate}
\end{lemma}

Let 
$$V:=\left\{q\in (1, M+1]:\overline{\alpha(q)} \preceq\alpha_{i+1}\alpha_{i+2}\cdots\preceq \alpha(q) \text{ for all } i\geq 0\right\}.$$ The set $V$ has zero Lebesgue measure, see \cite{DK2009}. 
For each $M\in \mathbb{N}$ we define the generalised golden ratio $q_{GR}$ to be the unique $q\in (1,M+1)$ for which $(0,M/(q-1))\cap \uu_{q}\neq \emptyset$ for $q>q_{GR}$ and $(0,M/(q-1))\cap \uu_{q}= \emptyset$ for $q<q_{GR}$. In \cite{B2014} it was shown that for each $M\in\mathbb{N}$ a generalised golden ratio exists, and is given by the following formula:
\[ q_{GR} = \left\{ \begin{array}{ll}
k+1 & \mbox{if $M=2k$};\\
\frac{k+1+\sqrt{k^2+6k+5}}{2} & \mbox{if $M=2k+1$}.\end{array} \right. \] $q_{GR}$ is the smallest element of $V$ \cite{DKL2016} and $M+1$ is  the
largest elements of $V$. By \cite[Theorem 1.3]{DKL2016}, we have 
$$[q_{GR}, M+1]\setminus V=(q_{GR}, M+1)\setminus V=\bigcup(q_\ell, q_r),$$
where the union on the right hand-side is pairwise disjoint and countable.  The open
intervals $(q_\ell, q_r)$ are referred to as the basic intervals of $[q_{GR}, M+1]\setminus V$. The following property of basic intervals was established in \cite{DK2009}.

\begin{lemma}\label{l22qqq}
Let $(q_l,q_r)$ be a basic interval. For any $q_1,q_2 \in (q_l,q_r)$ we have $U_{q_1}=U_{q_2}$.
\end{lemma}

We now recall some technical results on $\uu_{q}^{j}$. The following lemma follows from \cite[Lemma 1.6]{BS2014}. We remark that this lemma is phrased in the case when $M=1$ but the same argument applies for arbitrary $M$.
\begin{lemma}\label{l22qq}
Let $q>1$ and $x\in \uu_q^j$ for some $j\geq 3$.
Then there exists $(c_i)$ a $q$-expansion of $x$ and an integer $k\geq 0$ such that $\pi_q((c_{k+i}))\in \uu_q^2$.
\end{lemma}
Lemma \ref{l22qq} implies the following statement.
\begin{lemma}\label{l53}
For any $q>1$ and $j\geq 3,$ we have
\begin{equation*}
\dim_H\uu_q^j\le\dim_H\uu_q^2.
\end{equation*}
\end{lemma}

The following lemma follows from results proved in \cite{B2014}. See the discussion after Lemma 2.8 in this paper.

\begin{lemma}\label{l22q}
Assume $x=\pi_q((c_i))\in \uu_q\subset I_q$.
\begin{enumerate}[\upshape (i)]
\item $c_1=0$ if and only if $x\in[0,1/q)$.
\item For $i=1,\ldots, M-1,$ $c_1=i$ if and only if $x\in \left(\frac{(i-1)(q-1)+M}{q^2-q}, \frac{i+1}{q}\right).$
\item $c_1=M$ if and only if $x\in\left(\frac{(M-1)(q-1)+M}{q^2-q}, \frac{M}{q-1}\right]$ .
\end{enumerate}
\end{lemma}
With Lemma \ref{l22q} in mind, we define the switch region as follows.
$$S_q:=\bigcup_{i=1}^M\left[\frac{i}{q}, \frac{(i-1)(q-1)+M}{q^2-q}\right].$$ Note that $S_{q}$ is the complement to the intervals listed in items (i), (ii) and (iii) in Lemma \ref{l22q}. The following properties of the switch region were established in \cite{B2014}.
\begin{lemma}\label{switch region lemma}
Let $q\in (q_{GR},M+1]$. Then the following statements are true:
\begin{enumerate}[\upshape (i)]
	\item Let $x\in I_{q}$. Then there exists $(a_i),(b_i)\in \{0,\ldots,M\}^{\mathbb{N}}$ such that  $\pi_q((a_i))=\pi_q((b_i))=x$ and $a_1\neq b_1$ if and only if $x\in S_{q}$.
	\item Let $k\in \{1,\ldots,M\}$ and $x\in I_{q}.$ Then there exists $(a_i),(b_i)\in \{0,\ldots,M\}^{\mathbb{N}}$ such that $a_1=k$, $b_1=k-1$ and $\pi_q((a_i))=\pi_q((b_i))=x$ if and only if $x\in \left[\frac{k}{q}, \frac{(k-1)(q-1)+M}{q^2-q}\right].$
	\item If $i\neq j$ then $\left[\frac{i}{q}, \frac{(i-1)(q-1)+M}{q^2-q}\right]\cap \left[\frac{j}{q}, \frac{(j-1)(q-1)+M}{q^2-q}\right]=\emptyset.$
\end{enumerate}
\end{lemma}
Lemma \ref{switch region lemma} implies the useful fact that if $x\in I_{q}$ and there exists $(a_i),(b_i)\in \{0,\ldots,M\}^{\mathbb{N}}$ such that $\pi_q((a_i))=\pi_q((b_i))=x$ and $a_1\neq b_1,$ then $a_1$ and $b_1$ are successive digits in $\{0,\ldots,M\}$, and if $(c_i)$ is another sequence such that $\pi_q((c_i))=x$ then either $c_1=a_1$ or $c_1=b_1$.
\section{Properties of $\uu_q^2$}\label{s3}
In this section we prove some properties of the set $\uu_q^2$ and associated power series. 
\begin{lemma}\label{l31qqq}
Let $q\in (q_{GR}, M+1]$. Then $x\in \uu_q^2$ if and only if there exist $(a_i), (b_i)\in U_{q}$, a finite word $w\in \cup_{n=0}^{\infty}\{0,1,\ldots, M\}^n$ and $0\leq m<M$
 such that
 \begin{equation}\label{e31qq}
 x=\pi_q(wm(a_i))=\pi_q(w(m+1)(b_i))
 \end{equation}
 and \begin{equation}\label{e32qq}
 \pi_{q}(w_{j+1}\ldots w_nm(a_i))\notin S_{q}\text{ for all }  0\leq j< n.
 \end{equation}
\end{lemma}
\begin{proof}
Let $x\in I_{q}$ and suppose that \eqref{e31qq} and \eqref{e32qq} hold for some $(a_i), (b_i)\in U_{q}$. Then \eqref{e31qq} implies that
 $x$ has at least two different $q$-expansions, which are $wm(a_i)$ and $w(m+1)(b_i)$. Furthermore, Lemma \ref{switch region lemma}, \eqref{e32qq}, and the fact $(a_i), (b_i)\in U_{q}$ imply that these are the only $q$-expansions. Therefore $x\in \uu_{q}^2$. 

Now we prove the necessity. Take $x\in \uu_q^2$ having exactly two $q$-expansions $(c_i)$ and $(d_i)$. Let $n\geq1$ be the least integer such that $c_1\cdots c_{n-1}=d_1\cdots d_{n-1}$ and $c_n\neq d_n$. Without loss of generality we assume $c_n<d_n$. By Lemma \ref{switch region lemma} we know that there exists $0\leq m<M$ such that 
$$\pi_q((c_{n-1+i}))=\pi_q((d_{n-1+i}))\in \left[\frac{m+1}{q},\frac{m(q-1)+ M}{q^2-q}\right]\textrm{ and }c_n=m,\, d_n=m+1. $$
Moreover, because $x\in \uu_{q}^{2}$ we know that $(c_{n+i}), (d_{n+i})\in
U_q.$ Taking $w=c_1\cdots c_{n-1},$ $(a_i)=(c_{n+i})$, and $(b_i)=(d_{n+i}),$ we see that $ x=\pi_q(wm(a_i))=\pi_q(w(m+1)(b_i))$ and so \eqref{e31qq} holds. To see that \eqref{e32qq} holds, we remark that if it did not hold then Lemma \ref{switch region lemma} would imply that we have a choice of digit before the $n$-th position. This would imply that $x$ has at least three different $q$-expansions and so would contradict that $x\in \uu_{q}^{2}$. Therefore \eqref{e32qq} holds.
\end{proof}

 As we will see later on in the paper, two sequences $(a_i), (b_i)\in U_{q}$ give rise to a point with exactly two $q$-expansions if and only if $q$ is the zero of some appropriate power series with coefficients given by $(a_i)$ and $(b_i)$. We now set out to show that such a power series has at most one zero in $(q_{KL},M+1)$. To prove this we will make use of ideas from \cite{PerSol}.\\

Fix $M\geq 1.$ We consider functions of the form 
\begin{equation}
\label{gform}
g(x)=1+\sum_{i=1}^{\infty}b_ix^i\text{ with } b_i\in \{0, \pm 1, \cdots \pm M \}.
\end{equation} We say that the \emph{$\delta$-transversality condition} holds on the interval $I$ for some $\delta>0$ if for any $x\in I$ and $g$ of the form \eqref{gform}, whenever $g(x)<\delta$ we have $g'(x)<-\delta$. A power series $h$ is called  a \emph{$(*)$-function} if for some $k\geq 1$ and $a_k\in[-M, M]$ we have
\begin{equation*}
h(x)=1-M\sum_{i=1}^{k-1}x^i+a_k x^k+M\sum_{i=k+1}^\infty x^i.
\end{equation*} The following lemma connects $(*)$-functions and $\delta$-transversality.

  \begin{lemma}\label{l13}
	\label{*-function lemma}
If $h$ is a $(*)$-function such that
\begin{equation*}
h(x_{M})>\delta\;\; \text{and}\;\; h'(x_{M})<-\delta
\end{equation*}
for some $x_{M}\in (0,1)$ and $\delta\in(0,1)$, then the $\delta$-transversality condition holds on the interval $[0,x_{M}]$.
\end{lemma}
\begin{proof}
The case where $M=1$ was proved in \cite{PerSol}. The case where $M>1$ is proved by an analogous argument. 

\end{proof}
The following lemma will allow us to establish $\delta$-transversality within the interval $[0,q_{KL}^{-1}]$ for some $\delta>0$.



\begin{lemma}\label{*-function lemma 2}
For each $M\in\mathbb{N}$ let the ($*$)-function $h$ and $x_{M}$ be defined as follows:
\begin{itemize}
	\item Assume $M$ is of the form $M=2k+1$. Let
	\begin{equation*}
	h(x)=
	\begin{cases}1-x-x^2-x^3+0.5x^4+\sum_{i=5}^\infty x_i&\text{when $k=0$,}\\
	1-3x-0.5x^2+3\sum_{i=3}^\infty x^i&\text{when $k=1$,}\\
	1-(k+3)x+(2k+1)\sum_{i=2}^\infty x^i &\text{when $k\geq 2$.}\\
	\end{cases}
	\end{equation*}Let
	\begin{equation*}
	x_{M}=
	\begin{cases}2^{-2/3}&\text{when $k=0$,}\\
	\frac{-(k+1)+\sqrt{(k+1)^2+4(k+1)}}{2(k+1)}&\text{when $1\leq k<3$,}\\
	1/(k+1) &\text{when $k\geq 3$.}\\
	\end{cases}
	\end{equation*}
	\item Assume $M$ is of the form $M=2k$. Let
	\begin{equation*}
	h(x)=
	\begin{cases}
	1-2x-0.5x^2+2\sum_{i=3}^\infty x^i&\text{whenever $k=1$,}\\
	1-(k+2)x+2k\sum_{i=2}^\infty x^i &\text{whenever $k\geq 2$.}\\
	\end{cases}
	\end{equation*}
	Let
	\begin{equation*}
	x_{M}=
	\begin{cases}	\frac{-(k+1)+\sqrt{(k+1)^2+4k}}{2k}&\text{when $1\leq k<3$,}\\
	1/(k+1) &\text{when $k\geq 3$.}\\
	\end{cases}
	\end{equation*}
\end{itemize} Then $h(x_{M})>0$ and $h'(x_{M})<0$ and  $[1/(M+1), 1/q_{KL}]\subset [0, x_M].$
\end{lemma}

\begin{proof}
We begin by remarking that the fact $h(x_{M})>0$ and $h'(x_{M})<0$ was established when $M=1$ in \cite{PerSol}. Moreover, $[1/2, 1/q_{KL}]\subset [0, 2^{-2/3}]$ follows from the fact $q_{KL}\approx 1.787\ldots$. 

We now prove that $x_M^{-1}<q_{KL}$ for all  $M\geq 2$. This implies our final assertion $[1/(M+1), 1/q_{KL}]\subset [0, x_M]$ for the remaining values of $M$. When $M=3$ or $5$, it follows from the definition that $x_{M}$ is the solution to $(k+1)x^2+(k+1)x-1=0$ for the appropriate value of $k$. Using this algebraic relation it can be shown that the quasi-greedy expansion of $1$ in base $x_{M}^{-1}$ is $((k+1)k)^\infty.$ Similarly, it can be shown when $M=2$ or $M=4$ that the quasi-greedy expansion of $1$ in base $x_{M}^{-1}$ is $((k+1)(k-1))^\infty$ for the appropriate value of $k$. For the remaining vales of $M$ the quasi-greedy expansion of $1$ in base $x_{M}^{-1}$ is $k^{\infty}$. By inspection of \eqref{12} we see that the quasi-greedy expansion of $1$ in base $x_M^{-1}$ is lexicographically smaller than $\alpha(q_{KL})$ for any $M\in \mathbb{N}$. Lemma \ref{quasi-greedy lemma} (i) then implies $x_M^{-1}<q_{KL}$.




It remains to show that $h(x_M)>0$ and $h'(x_M)<0$. For $M=2,3,4,5$ this can be checked by a direct computation. For $M=2k+1$ and $k\ge 3$,  we have
\begin{equation*}
h'(x)=-3-k+(2k+1)\frac{2x-x^2}{(x-1)^2}.
\end{equation*}
Substituting in our value for $x_{M}$ we have
\begin{equation*}
h(x_M)=\frac{1}{k+k^2}\qtq{and}
h'(x_M)=1+\frac{1}{k^2}+\frac{4}{k}-k.
\end{equation*}
We always have $h(x_M)>0$ and it is simple to show that  $h'(x_M)<0$ for all $k\geq 3$.

For $M=2k$ and $k\ge 3$ we have
\begin{equation*}
h'(x)=-2-k+2k\frac{2x-x^2}{(x-1)^2}.
\end{equation*}
Substituting in our value of $x_{M}$ we have
\begin{equation*}
h(x_M)=\frac{1}{k+1}\qtq{and}
h'(x_M)=2+\frac{2}{k}-k.
\end{equation*}
We always $h(x_M)>0$ and it is easy to show that  $h'(x_M)<0$ for all $k\geq 3$.
\end{proof}

\begin{proposition}\label{p35}
Let $(a_i), (b_i)\in \{0, 1,\cdots, M\}^\mathbb{N}$. Then there exists at most one $q\in [q_{KL}, M+1]$ such that $\pi_{q}((a_i))=\pi_{q}((b_i))+1$.
\end{proposition}
\begin{proof}
 The equation $\pi_{q}((a_i))=\pi_{q}((b_i))+1$ can be rewritten as $0=1+\sum_{i=1}^{\infty}\frac{b_i-a_i}{q^i}.$ The zeros of this function can be mapped to the zeros of $0=1+\sum_{i=1}^{\infty}(b_i-a_i)x^i$ by taking this reciprocal. Importantly this map is a bijection. By Lemmas \ref{*-function lemma} and \ref{*-function lemma 2} we know that the $\delta$-transversality condition holds on the interval $[1/(M+1),q_{KL}^{-1}]$ for $\delta=\min\{h(x_M)/2,|h'(x_M)|/2\}.$ Therefore there exists at most one $x^*\in [1/(M+1),q_{KL}^{-1}]$ such that $0=1+\sum_{i=1}^{\infty}(b_i-a_i)(x^*)^i.$ Using the bijection mentioned above we may conclude our result.
\end{proof}



\section{Proofs of Theorems \ref{t11} and \ref{t12}}\label{s4}
The main focus of this section will be to prove Theorem \ref{t11}. We will then explain how the argument can be adapted to prove Theorem \ref{t12}. We start by remarking that by Lemma \ref{l53} we know that for any $q>1$ we have $$\dim_H\uu_q^j\leq \dim_H\uu_q^2$$ for every $j\in\{3,4,\ldots\}$.
Because of this lemma, to prove Theorem \ref{t11} it suffices to show that for Lebesgue almost every $q\in (q_{KL},M+1)$ we have $$\dim_{H}\uu_{q}^{2}\leq \max\{2\dim_{H}\uu_{q}-1,0\}.$$ Moreover, to prove this statement it is sufficient to show that for any $\epsilon>0$ \begin{equation}
	\label{epsilon bound}	
	\dim_{H}\uu_{q}^{2}\leq \max\{2\dim_{H}\uu_{q}-1+\epsilon,0\}
	\end{equation} 
	for Lebesgue almost every $q\in (q_{KL},M+1)$. We now fix such an $\epsilon>0$ and set out to show that this is the case.\\

To prove \eqref{epsilon bound} holds for Lebesgue almost every $q\in (q_{KL},M+1)$ we will make use of a countable collection of closed intervals $\{J_i\}$ that are chosen in a way that depends upon $\epsilon$ and $M$. These intervals satisfy the following properties:

\begin{enumerate}[\upshape (i)]
\item $\{J_i\}$ cover $(q_{KL}, M+1)$ up to a set of Lebesgue measure zero.
\item Each $J_i$ is small in a way that depends upon $\epsilon$ and $M$. 
\item Each $J_i$ is contained in a basic interval $(q_\ell, q_r)$ of $[q_{GR}, M+1]\setminus V$. 
\item When $M\geq 2$, for each $J_i$ there exists $p\in\mathbb{N}$ such that $J_i\subset (p,p+1).$
\item For any $J_i$, if $q_1,q_2\in J_i$ then $|\dim_{H}\uu_{q_1}-\dim_{H}\uu_{q_2}|<\epsilon/2.$
\end{enumerate}

Item (ii) is important because we will soon show that if the $\{J_i\}$ are chosen to be small in an appropriate way that depends upon $\epsilon$ and $M$, then certain useful properties hold. We remark that item (iii) in the above together with Lemma \ref{l22qqq} implies that for any $q_1, q_2 \in J_i\subset (p_\ell, p_r)$ we have $$U_{q_1}=U_{q_2}.$$ Item (iv) is a technical assumption that will help with our later proof. Item (v) follows from the continuity of the function which maps $q$ to $\dim_{H}\uu_q.$ This fact was established in \cite{KKL2017,AK2019}.
\color{black}


From the collection $\{J_i\}$ we now make an arbitrary choice that we will denote by $J$. For the rest of the proof of Theorem \ref{t11} the interval $J$ is fixed. By item (i) from the properties of $\{J_i\}$ mentioned above, and the arbitrariness of $J,$ to prove \eqref{epsilon bound} holds for Lebesgue almost every $q\in (q_{KL},M+1)$ it is sufficient to show that it holds for Lebesgue almost every $q\in J$. To do this, we need to introduce two new sets and a function between them. 




Let $q_{0}:=\min\{q:q\in J\}.$ To each $0\leq m<M$ we associate $$A_{m}:=\bigcup_{\stackrel{a, b\in U_{ q_0}}{\exists q^*\in J: \pi_{q^*}(a)=\pi_{q^*}(b)+1}}\left(\pi_{q_0}(ma),\pi_{q_0}((m+1)b\right)$$ and
\begin{equation}\label{A_q}
A:=\bigcup_{m=0}^{M-1}A_{m}.
\end{equation} 
We emphasise that for each $m$ the terms in the union in $A_{m}$ are elements of $\mathbb{R}^2$ and not intervals. Even though $A$ consists of elements of $\mathbb{R}^2$ obtained by applying $\pi_{q_0}$ to a collection of sequences from $U_{q_0}$, because $U_{q}=U_{q_0}$ for each $q\in J$ the set $A$ in fact contains information about the whole interval $J$ and not just the specific base $q_0$. This fact is crucial in our proof. By Proposition \ref{p35}, we know that if $a,b\in U_{q_0}$ are such that there exists $q^*\in J$ for which $\pi_{q^*}(a)=\pi_{q^*}(b)+1,$ then $q^*$ is unique. When such a $q^*$ exists then we will denote it by $q^*(a,b).$ 

To each $0\leq m<M$ we associate
$$B_{m}:=\bigcup_{\stackrel{a,b\in U_{q_0}}{\exists q^*\in J: \pi_{q^*}(a)=\pi_{q^*}(b)+1}}(q^{*}(a,b),\pi_{q^*}(ma))$$ and 
\begin{equation*}
B:=\bigcup_{m=0}^{M-1}B_{m}.
\end{equation*} 
As above, for each $m$ the terms in the union in $B_{m}$ are elements of $\mathbb{R}^2$ and are not intervals. We emphasise that both $A$ and $B$ depend upon our choice of interval $J$. However because $J$ is fixed we suppress this dependence within our notation. 

The following lemma connects vertical slices through $B$ with the set $\uu_{q}^2.$
\begin{lemma}\label{l43}
	For any $q\in J$ we have $$\dim_H\left(\{(q,y):y\in\mathbb{R}\}\cap B \right)=\dim_H\uu_{q}^2.$$
	Moreover, if $\{(q,y):y\in\mathbb{R}\}\cap B=\emptyset$, then $\uu_{q}^2=\emptyset$.
\end{lemma}
\begin{proof}
	Fix $q\in J$. We begin by remarking that 
	$$
	\left\{(q,y): y\in\mathbb{R}\right\}\cap B=\bigcup_{m=0}^{M-1}\bigcup_{\stackrel{a,b\in U_{q_0}}{ \pi_{q}(a)=\pi_{q}(b)+1}}(q,\pi_{q}(ma)).
	$$ Moreover, because $q\in J$ and $U_{q_0}=U_{q}$ for all $q\in J$ we in fact have $$
	\{(q,y): y\in\mathbb{R}\}\cap B=\bigcup_{m=0}^{M-1}\bigcup_{\stackrel{a,b\in U_{q}}{ \pi_{q}(a)=\pi_{q}(b)+1}}(q,\pi_{q}(ma)).
	$$
	Therefore 
	\begin{equation}
	\label{dimension equality}
	\dim_H\left(\{(q,y): y\in \mathbb{R}\}\cap B\right)=\dim_H C_{q}
	\end{equation}
	where $$C_{q}:=\bigcup_{m=0}^{M-1}\bigcup_{\stackrel{a,b\in U_{q}}{ \pi_{q}(a)=\pi_{q}(b)+1}}\pi_{q}(ma).$$
	If $a,b\in U_q$ are such that $\pi_{q}(a)=\pi_{q}(b)+1$ then $\pi_{q}(ma)=\pi_{q}((m+1)b)$ for each $0\leq m<M.$ Therefore by Lemma \ref{l31qqq} we may conclude that
	$\pi_{q}(ma)\in \uu_{q}^2$ for all $0\leq m<M$ and $$\dim_H C_{q}\leq\dim_H\uu_{q}^2.$$

	Now we prove that $\dim_H C_{q}\geq\dim_H\uu_{q}^2$. Suppose $x\in \uu_{q}^2$, then by Lemma \ref{l31qqq} again, there exists a finite word $w\in \cup_{n=0}^{\infty} \{0,\cdots, M\}^n$, $0\leq m<M$ and $a, b\in U_{q}$ such that  $x=\pi_{q}(wma)=\pi_{q}(w(m+1)b)$.
	Therefore
	\begin{equation*}
	x\in \pi_{q}(w0^\infty)+\frac{1}{q^n} C_{q},
	\end{equation*}
	and so
	\begin{equation*}
	\uu_{q}^2\subseteq\bigcup_{n=0}^{\infty}
	\bigcup_{w\in\set{0,\ldots,M}^n}
	\left(\pi_{q}(w0^{\infty})+\frac{1}{q^n} C_{q}\right).
	\end{equation*}
	Since $\uu_{q}^2$ is covered by countably many sets, each of which is an affine image of $ C_{q}$,  we conclude that $\dim_H\uu_{q}^2\le\dim_H  C_{q}$. We have shown that $\dim_H\uu_{q}^2=\dim_H  C_{q}$. Combining this equality with \eqref{dimension equality}  we may conclude that $\dim_H\left(\{(q,y):y\in\mathbb{R}\}\cap B \right)=\dim_H\uu_{q}^2$ as required.
	
	It remains to verify that if $\{(q,y):y\in\mathbb{R}\}\cap B=\emptyset$ then $\uu_{q}^2=\emptyset$. If $\{(q,y):y\in\mathbb{R}\}\cap B=\emptyset$ then $C_{q}$ as defined above is also empty. The fact $\uu_{q}^2$ is covered by countably many affine images of $C_{q}$ holds even if $C_{q}$ is empty. Therefore if $\{(q,y):y\in\mathbb{R}\}\cap B$ is empty then so is $\uu_{q}^2$.
	
\end{proof}

What remains of our proof of Theorem \ref{t11} focuses on the properties of a map $f:A\to B$ defined below. To make sure this map is well defined it is necessary to assume that $A_{i}\cap A_{j}=\emptyset$ for $i\neq j$. The following lemma allows us to make such an assumption.
\begin{lemma}\label{l41q}
The intervals $\{J_i\}$ can be chosen to be sufficiently small such that $A_{i}\cap A_{j}=\emptyset$ for $i\neq j$.
\end{lemma}
\begin{proof}
Consider the interval $$\left[\frac{i+1}{q_0},\frac{i( q_0-1)+M}{q_0^2- q_0}\right].$$ If $i\neq j$ then Lemma \ref{switch region lemma} tells us that
$$\left[\frac{i+1}{ q_0},\frac{i(q_0-1)+M}{ q_0^2- q_0}\right]\cap \left[\frac{j+1}{ q_0},\frac{j( q_0-1)+M}{ q_0^2- q_0}\right]=\emptyset.$$ Moreover, there exists $\delta>0$ (for example we could take $\delta:= \frac{1}{3}\frac{2q_{KL}-2-M}{q_{KL}^2-q_{KL}}$) depending only upon $M$ such that 
\begin{equation}
\label{intersection}
\left[\frac{i+1}{q_0}-\delta,\frac{i(q_0-1)+M}{q_0^2-q_0}+\delta\right]\cap \left[\frac{j+1}{q_0}-\delta,\frac{j( q_0-1)+M}{ q_0^2-q_0}+\delta\right]=\emptyset.
\end{equation} 
Let $a,b\in U_{q_0}$ be such that there exists $q^{*}\in J$ such that $\pi_{q^*}(a)=\pi_{q^*}(b)+1.$ Then for each $0\leq i<M$ we have $\pi_{q^*}(ia)=\pi_{q^*}((i+1)b)$ and $$\pi_{q^*}(ia)\in \left[\frac{i+1}{q^{*}},\frac{i(q^{*}-1)+M}{(q^*)^2-q^*}\right]$$ by Lemma \ref{switch region lemma}. By continuity, we can assume that our intervals $\{J_i\}$ were chosen to be sufficiently small such that 
\begin{equation*}
\left[\frac{i+1}{q^{*}},\frac{i(q^{*}-1)+M}{(q^*)^2-q^*}\right]\subset \left(\frac{i+1}{q_0}-\delta,\frac{i(q-1)+M}{ q_0^2- q_0}+\delta\right)
\end{equation*}for any $q^{*}\in J$. Therefore, by continuity we can choose our $\{J_i\}$ to be sufficiently small such that
\begin{equation}
\label{inclusion}
\pi_{q_0}(ia),\pi_{q_0}((i+1)b)\in\left[\frac{i+1}{q_0}-\delta,\frac{i(q_0-1)+M}{ q_0^2- q_0}+\delta\right].
\end{equation} Combining \eqref{intersection} and \eqref{inclusion} we may conclude that $A_{i}\cap A_{j}$ for $i\neq j$.
\end{proof}

For each $0\leq m<M$ we define $f_{m}:A_{m}\to B_{m}$ by $$f_{m}\left((\pi_{q_0}(ma),\pi_{q_0}((m+1)b))\right)=(q^{*}(a, b),\pi_{q^*}(ma)).$$ Since we can assume the $A_{m}$ are disjoint by Lemma \ref{l41q}, we can define \begin{equation}\label{f function}
f:A\to B
\end{equation}
by $f|_{A_{m}}=f_{m}.$ Importantly, it follows from the definitions of $A$, $B$ and $f$ that $f(A)=B$. Our goal now is to prove that the function $f:A\to B$ is Lipschitz. The following three lemmas are all technical results that will allow us to establish this fact. 

\color{black}
\begin{lemma}\label{l41qq}
Let $q\in(q_{KL}, M+1)$ and  $(a_i), (b_i)\in U_{q}$. If
\begin{equation}\label{e41qqq}\pi_{q}((a_i))=\pi_{q}((b_i)) +1\end{equation} then for all $n\geq 1$ we have
\begin{enumerate}[\upshape (i)]
\item 
$
(a_{i})\succeq\alpha(q)\succ \alpha(q_{KL})\text{ and }
\pi_q(a_1\cdots a_n 0^\infty)\geq\pi_q(\alpha_1(q_{KL})\cdots \alpha_n(q_{KL})0^\infty).
$

\medskip
\item $
(b_{i})\preceq\overline{\alpha(q)}\prec \overline{\alpha(q_{KL})}
\text{ and }
\pi_q(b_1\cdots b_n M^\infty)\leq\pi_q(\overline{\alpha_1(q_{KL})\cdots \alpha_n(q_{KL})}M^\infty).
$

\end{enumerate}
\end{lemma}
\color{black}
\begin{proof}
Let $(a_i), (b_i)\in U_{q}$ satisfy \eqref{e41qqq}. Equation \eqref{e41qqq} implies that $\pi_{q}((a_i))\geq 1.$ Since $(a_i)\in U_{q}$ it must be the greedy expansion of some element in $I_q$ greater than or equal to $1$. By Lemma \ref{quasi-greedy lemma} (ii) we have $$(a_i)\succeq\alpha(q).$$ Moreover, because $q\in (q_{KL},M+1)$, we know by Lemma  \ref{quasi-greedy lemma} (i) that $\alpha(q)\succ \alpha(q_{KL}).$ This implies the first part of (i). 

Now we show the second part of (i). We  claim that the two sequences $a_1\cdots a_n 0^\infty$ and $\alpha_1(q_{KL})\cdots \alpha_n(q_{KL})0^\infty$ are greedy $q$-expansions. This is true because of the characterisation of greedy expansions in base $q$ provided by Lemma \ref{quasi-greedy lemma} (ii), and because both $(a_i)$ and $\alpha(q_{KL})$ are greedy expansions in base $q$. It now follows from 
$(a_{i})\succ\alpha(q_{KL})$ that $$a_1\cdots a_n 0^\infty \succeq \alpha_1(q_{KL})\cdots \alpha_n(q_{KL})0^\infty.$$
Thus the second part of (i) follows from Lemma \ref{quasi-greedy lemma} (iii).

Statement (ii) is proved similarly, this time exploiting the fact that if $(a_i), (b_i)\in U_{q}$ satisfy \eqref{e41qqq} then we must have $\pi_{q}((b_i))\leq \frac{M}{q-1}-1.$

\end{proof}

\begin{lemma}\label{l45qq}
	We can  choose our $\{J_i\}$ in such a way so that the following is true: There exists $C_0>0$ that does not depend upon our choice $J$, such that if $(a_i), (b_i)\in U_{q_0}$ are such that 
\begin{equation}\label{e410qq}
	\pi_{q}((a_i))=\pi_{q}((b_i))+1
\end{equation} 
for some $q\in J$, then for all $n$ sufficiently large we have $$\inf_{q_1,q_2\in J}\left|\frac{1}{q_{1}q_{2}}+\sum_{i=1}^{n}\frac{(b_i-a_i)(q_{2}^{i}+q_{2}^{i-1}q_{1}+\cdots +q_{1}^i)}{q_{1}^{i+1}q_{2}^{i+1}}\right|\geq C_0. $$
	\end{lemma}
\begin{proof}
	If we set $q_{2}=q_{1}$ in the above expression we obtain $$\inf_{q_{1}\in J}\left|\frac{1}{q_{1}^{2}}+\sum_{i=1}^{n}\frac{(b_i-a_i)(i+1)}{q_{1}^{i+2}}\right|.$$
Now suppose that we have shown that there exists $C_1>0$ independent of 
$J,$ such that if $(a_i), (b_i)\in U_{q_0}$ are such that $$\pi_{q}((a_i))=\pi_{q}((b_i))+1$$ for some $q\in J,$ then for all $n$ sufficiently large we have
	\begin{equation}
	\label{need to show}
	\inf_{q_{1}\in J}\left|\frac{1}{q_{1}^{2}}+\sum_{i=1}^{n}\frac{(b_i-a_i)(i+1)}{q_{1}^{i+2}}\right|\geq C_1.
	\end{equation} Then by a continuity and compactness argument, if the  $\{J_i\}$ are chosen sufficiently small, then $$\inf_{q_1,q_2\in J}\left|\frac{1}{q_{1}q_{2}}+\sum_{i=1}^{n}\frac{(b_i-a_i)(q_{2}^{i}+q_{2}^{i-1}q_{1}+\cdots +q_{1}^i)}{q_{1}^{i+1}q_{2}^{i+1}}\right|\geq \frac{C_1}{2}$$ for all $n$ sufficiently large. Therefore to complete our proof we only need to show that inequality \eqref{need to show} is satisfied for all $n$ sufficiently large. \\

We can choose our $\{J_i\}$ to be sufficiently small so that for all $n$ sufficiently large, we have by \eqref{e410qq} that
\begin{equation}\label{e411qqqq}
\sum_{i=1}^n\frac{a_i}{q_1^i}-\sum_{i=1}^n\frac{b_i}{q_1^i}\geq 1-\omega.
\end{equation}
for all $q_1\in J$. Here we fix $\omega=0.01$.

Using \eqref{e411qqqq} for $n$ sufficiently large we have 
\begin{equation}\label{e412qqq}
\begin{split}
\frac{1}{q_{1}^{2}}+\sum_{i=1}^{n}\frac{(b_i-a_i)(i+1)}{q_{1}^{i+2}}&=\frac{1}{q_{1}^{2}}+2\sum_{i=1}^{n}\frac{b_i-a_i}{q_{1}^{i+2}}+\sum_{i=1}^{n}\frac{(b_i-a_i)(i-1)}{q_{1}^{i+2}}\\
&\leq \frac{1}{q_{1}^{2}}(-1+2\omega)+\sum_{i=1}^{n}\frac{(b_i-a_i)(i-1)}{q_{1}^{i+2}}
\end{split}
\end{equation}
We now proceed via a case analysis.

\noindent \textbf{Case 1. $M=1$. }
It follows from   Lemma \ref{l41qq} and \eqref{12} that 
\begin{equation}\label{e413qqq}
(a_i)\succ\alpha(q_{KL})=11010011\cdots\text{ and } (b_i)\prec \overline{\alpha(q_{KL})}=00101100\cdots,
\end{equation}
then we  obtain
 \begin{equation}\label{e46qq}b_1-a_1=-1, b_2-a_2=-1 \text{ and } b_i-a_i\leq 1\text{ for all } i\geq 3.
 \end{equation}
\noindent \textbf{Case 1a.} If  $b_3-a_3\leq 0$, then by \eqref{e412qqq} and \eqref{e46qq} we obtain
\begin{align*}
\frac{1}{q_{1}^{2}}+\sum_{i=1}^{n}\frac{(b_i-a_i)(i+1)}{q_{1}^{i+2}}
&\leq \frac{1}{q_{1}^{2}}(-1+2\omega)+\sum_{i=1}^{n}\frac{(b_i-a_i)(i-1)}{q_{1}^{i+2}}\\
&< \frac{1}{q_{1}^{2}}(-1+2\omega) -\frac{1}{q_1^4}+\frac{0}{q_1^5}+\sum_{i=4}^{\infty}\frac{i-1}{q_{1}^{i+2}}\\
&=\frac{1}{q_{1}^{2}}(-1+2\omega)-\frac{2}{q_1^4} -\frac{2}{q_1^5}+\frac{1}{q_1^2(q_1-1)^2}.\\
\end{align*}
The last equality follows from $$\frac{1}{(q_1-1)^2}=\sum_{i=1}^\infty\frac{i-1}{q_1^{i}}.$$
A quick computer inspection verifies that $$\frac{1}{q_{1}^{2}}(-1+2\omega)-\frac{2}{q_1^4} -\frac{2}{q_1^5}+\frac{1}{q_1^2(q_1-1)^2}<0$$ for all $q_1\in (q_{KL}, 2)$. Therefore $C_1$ exists in this case. \\

\noindent \textbf{Case 1b.} If $b_3-a_3=1$, that is $b_3=1, a_3=0$. By \eqref{e413qqq} we must have  $b_4=0$ and $a_4=1$. Then it follows from   \eqref{e46qq} again that
\begin{align*}
\frac{1}{q_{1}^{2}}+\sum_{i=1}^{n}\frac{(b_i-a_i)(i+1)}{q_{1}^{i+2}}
&\leq \frac{1}{q_{1}^{2}}(-1+2\omega)+\sum_{i=1}^{n}\frac{(b_i-a_i)(i-1)}{q_{1}^{i+2}}\\
&< \frac{1}{q_{1}^{2}}(-1+2\omega) -\frac{1}{q_1^4}+\frac{2}{q_1^5}-\frac{3}{q_1^6}+\sum_{i=5}^{\infty}\frac{i-1}{q_{1}^{i+2}}\\
&=\frac{1}{q_{1}^{2}}(-1+2\omega)-\frac{2}{q_1^4} -\frac{6}{q_1^6}+\frac{1}{q_1^2(q_1-1)^2}.\\
\end{align*}
A quick computer inspection verifies that $$\frac{1}{q_{1}^{2}}(-1+2\omega)-\frac{2}{q_1^4} -\frac{6}{q_1^6}+\frac{1}{q_1^2(q_1-1)^2}<0$$ for all $q_1\in (q_{KL}, 2)$.
Therefore $C_1$ exists in this case.\\

\noindent \textbf{Case 2. $M\geq 2$.}
Let $\alpha(q_1)=(\alpha_i)$. Item (iv) in our list of assumption for the collection of interval $\{J_i\}$ tells us that $q_1\in J\subset (p,p+1]$ for some $p\leq M$. This implies that $\alpha_1\leq p$.
\color{black}
 We will now prove that
\begin{equation}\label{e48}b_i\leq p\text { for all  }i\geq 1.\end{equation} If $p=M$ then \eqref{e48} is trivially true. We now prove it for $p<M$.  

Combining $\alpha_1\leq p$ and Lemma \ref{quasi-greedy lemma} (i), we have  $\alpha(q_1)\preceq p^\infty$. By Lemma \ref{l41qq} we have $(b_{i})\preceq \overline{\alpha(q_1)}.$ Since $\alpha(q_1)\succeq \overline{\alpha(q_1)}$ for $q\geq q_{KL}$ we have $(b_{i})\preceq \alpha(q_1)\preceq p^\infty.$ Therefore $b_{1}\leq p$ and $b_1<M$. Moreover, since $(b_i)$ is the greedy expansion of some real number, Lemma \ref{quasi-greedy lemma} (iii) implies that $$(b_{n+i})\prec \alpha(q_1)$$ for all $n\geq 1$. This implies that $b_i\leq p$ for all $i\geq 2$. This completes our proof of \eqref{e48}.   

  Using \eqref{e412qqq} and \eqref{e48}, for $n$ sufficiently large we have 
\begin{align*}
\frac{1}{q_{1}^{2}}+\sum_{i=1}^{n}\frac{(b_i-a_i)(i+1)}{q_{1}^{i+2}}&=\frac{1}{q_{1}^{2}}+2\sum_{i=1}^{n}\frac{b_i-a_i}{q_{1}^{i+2}}+\sum_{i=1}^{n}\frac{(b_i-a_i)(i-1)}{q_{1}^{i+2}}\\
&\leq \frac{1}{q_{1}^{2}}(-1+2\omega)+\sum_{i=1}^{n}\frac{(b_i-a_i)(i-1)}{q_{1}^{i+2}}\\
&\leq  \frac{1}{q_{1}^{2}}(-1+2\omega) +p\sum_{i=1}^{\infty}\frac{(i-1)}{q_{1}^{i+2}}\\
&=\frac{1}{q_{1}^{2}}\left(-1+2\omega+\frac{p}{(q_{1}-1)^2}\right).
\end{align*} 
Recall that $\omega=0.01$. To complete our proof, it suffices to show that there exists $C_1>0$ such that  $$-1+2\omega+\frac{p}{(q_{1}-1)^2}<-C_{1}.$$ Suppose $p\geq 3.$ We must have $q_{1}>p$ since $J\subset (p,p+1)$. Therefore $$-1+2\omega+\frac{p}{(q_{1}-1)^2}\leq -1 +2\omega +\frac{p}{(p-1)^2}.$$ For $p\geq 3$ the right hand side of this expression can be uniformly bounded from above by a negative number. Therefore $C_{1}$ exists in this case. If $p=2$ then we are either in the case where $M=2$ or $M=3$. It can be verified using a computer that $$-1+2\omega+\frac{2}{(q_{1}-1)^2}<0$$ for any $q_{1}\geq 2.43$ for $M=2$ or $M=3$. One can use \eqref{e33qq} to verify that $q_{KL}>2.43$. The existence of $C_{1}$ follows.

\end{proof}

\begin{lemma}\label{l22}
There exists $C>0$ depending only upon J, such that for any $q\in J$ if $x=\pi_q((a_i)),y=\pi_q((b_i))\in \uu_q$ satisfies $$|x-y|\leq Cq^{-n},$$ then $b_i=a_i$ for all $1\leq i\leq n$.
\end{lemma}
\begin{proof}
Let $$C:=\min_{q\in J}\frac{1}{2}\left(\frac{M}{q-1}-1\right).$$
\color{black} We proceed by induction on $n$. For $n=1$ suppose $|x-y|\leq Cq^{-1},$ we now prove $b_1=a_1$. If $b_1\neq a_1$ then it follows from Lemma \ref{l22q} that $$|x-y|> \frac{M}{q^2-q}-\frac{1}{q}>Cq^{-1}.$$
Which is a contradiction.

Assume that if $|x-y|\leq Cq^{-n+1}$ then $b_i=a_i$ for all $1\leq i\leq n-1$. Using this assumption we now prove that if $|x-y|\leq Cq^{-n}$ then $b_i=a_i$ for all $1\leq i\leq n$.
Since $|x-y|\leq Cq^{-n}< Cq^{-n+1},$ our assumption on $n-1$ implies that  $b_i=a_i$ for all $1\leq i\leq n-1$.
Now we prove $b_n=a_n$. Assume that $b_n\neq a_n$. Then
\begin{equation}
\label{contradiction}
Cq^{-1}\geq q^{n-1}|x-y|=q^{n-1}\left|\sum_{i=1}^{\infty}\frac{a_i}{q^i}-\sum_{i=1}^{\infty}\frac{b_i}{q^i}\right|=\left|\sum_{i=1}^{\infty}\frac{a_{n-1+i}}{q^i}-\sum_{i=1}^{\infty}\frac{b_{n-1+i}}{q^i}\right|.
\end{equation} 
Because $a_{n}\neq b_n$, Lemma \ref{l22q} implies that $$\left|\sum_{i=1}^{\infty}\frac{a_{n-1+i}}{q^i}-\sum_{i=1}^{\infty}\frac{b_{n-1+i}}{q^i}\right|>Cq^{-1}.$$ This contradicts \eqref{contradiction} and so completes our proof.

\end{proof}

Equipped with Lemmas \ref{l45qq} and \ref{l22} we are now in a position to prove that $f$ is Lipschitz. To prove this statement it is convenient to use the infinity norm on $\mathbb{R}^2$ which we denote by $\|\cdot\|_{\infty}$ (i.e. $\|(x,y)\|_{\infty}:=\max\{|x|,|y|\}$). 

\begin{proposition}
    \label{p24}
Let $A$ and $f$ be defined as in  \eqref{A_q} and \eqref{f function}. Then there exists $C'>0$ depending only upon $J$  such that $$\|f(x_1,y_1)-f(x_2,y_2)\|_{\infty}\leq C'\|(x_1,y_1)-(x_2,y_2)\|_{\infty}$$ for all $(x_1,y_1),(x_2,y_2)\in A.$ 
\end{proposition}
\begin{proof}
To prove that such a $C'$ exists it suffices to consider $(x_1,y_1),(x_2,y_2)\in A$ for which $\|(x_1,y_1)-(x_2,y_2)\|_{\infty}$ is small. As such, by Lemma \ref{l41q} we can restrict our attention to those $(x_1,y_1),(x_2,y_2)\in A$ for which there exists a unique $m$ satisfying
$$(x_1,y_1),(x_2,y_2)\in \frac{\uu_{q_0}+m}{q_0}\times \frac{\uu_{q_0}+m+1}{q_0}.$$
Moreover by the definition of $A$ there exists $(a_i^1),(a_i^2),(b_i^1), (b_i^2)\in U_{q_0}$ such that
$$x_1=\pi_{q_0}(m(a_i^{1})), \quad x_2=\pi_{q_0}(m(a_i^{2})), \quad y_1=\pi_{q_0}((m+1)(b_i^{1})), \quad y_2=\pi_{q_0}((m+1)(b_i^{2})),$$  and $q_1,q_2\in J$ for which
\begin{equation}\label{e418}
\pi_{q_1}(m(a_i^{1}))=\pi_{q_1}((m+1)(b_i^{1}))\text{ and } \pi_{q_2}(m(a_i^{2}))=\pi_{q_2}((m+1)(b_i^{2})).
\end{equation}
By the definition of $f$, we have 
$$f((x_1, y_1))=(q_1, \pi_{q_1}(m(a_i^1)))\text{ and }f((x_2, y_2))=(q_2, \pi_{q_2}(m(a_i^2))).$$
Therefore  
\begin{equation*}
\begin{split}
\|f(x_1,y_1)-f(x_2,y_2)\|_{\infty}&=\|(q_1, \pi_{q_1}(m(a_i^1))-(q_2, \pi_{q_2}(m(a_i^2))\|_{\infty}\\
&=\max\left\{|q_1-q_2|, |\pi_{q_1}(m(a_i^1))-\pi_{q_2}(m(a_i^2))|\right\}.
\end{split}
\end{equation*}
We now bound the two terms appearing in this maximum accordingly. We begin by showing that there exists $C'$ such that $$|q_1-q_2|<C'\|(x_1,y_1)-(x_2,y_2)\|_{\infty}.$$
As stated above, it is sufficient to consider those $(x_1,y_1),(x_2,y_2)\in A$ for which $\|(x_1,y_1)-(x_2,y_2)\|_{\infty}$ is small. In particular, we can assume that $(x_1,y_1),(x_2,y_2)\in A$ are such that the unique $n$ satisfying
\begin{equation}
    \label{e26}
C q_0^{-n-2}<\|(x_1,y_1)-(x_2,y_2)\|_{\infty}\leq C q_0^{-n-1},
\end{equation}
is sufficiently large so that Lemma \ref{l45qq} applies. Here $C$ is as in Lemma \ref{l22}. Equation \eqref{e26} implies that $$|\pi_{q_0}((a_i^1))-\pi_{q_0}((a_i^2))|\leq C q_0^{-n} \textrm{ and }|\pi_{q_0}((b_i^1))-\pi_{q_0}((b_i^2))|\leq C q_0^{-n}.$$ Therefore Lemma \ref{l22} implies that $a_i^{1}=a_{i}^{2}$ and $b_i^{1}=b_{i}^2$ for $1\leq i\leq n.$

It follows from \eqref{e418} that $$\frac{1}{q_{1}}+\sum_{i=1}^{\infty}\frac{b_i^{1}-a_i^{1}}{q_{1}^{i+1}}=0\quad \textrm{ and }\quad \frac{1}{q_{2}}+\sum_{i=1}^{\infty}\frac{b_i^{2}-a_{i}^2}{q_{2}^{i+1}}=0.$$
Combining these two equations, we have
\begin{equation}\label{e21}
\frac{1}{q_{1}}+\sum_{i=1}^{\infty}\frac{b_i^{1}-a_i^{1}}{q_{1}^{i+1}}-\frac{1}{q_{2}}-\sum_{i=1}^{\infty}\frac{b_i^{2}-a_{i}^2}{q_{2}^{i+1}}=0.\end{equation} Using the fact $a_i^{1}=a_{i}^{2}$ and $b_i^{1}=b_{i}^2$ for $1\leq i\leq n,$ we see that \eqref{e21} yields
\begin{equation*}
    \begin{split}
&\left|\frac{1}{q_{1}}+\sum_{i=1}^{n}\frac{b_i^{1}-a_i^{1}}{q_{1}^{i+1}}-\frac{1}{q_{2}}-\sum_{i=1}^{n}\frac{b_i^{1}-a_{i}^1}{q_{2}^{i+1}}\right|\\
=&\left|\sum_{i=n+1}^{\infty}\frac{b_i^{2}
-a_{i}^2}{q_{2}^{i+1}}-\sum_{i=n+1}^{\infty}\frac{b_i^{1}-a_{i}^1}{q_{1}^{i+1}}\right|\\
\leq& \sum_{i=n+1}^{\infty}\frac{2M}{q_{0}^{i+1}}=\frac{2M}{q_0^2-q_0}q_0^{-n}<C_{1}q_{0}^{-n}.
    \end{split}
\end{equation*}
Here $C_1:=2M/(q_{KL}^2-q_{KL})$. In the final line we used that $q_0=\min\{q:q\in J\}$. It now follows from \eqref{e26} that 
\begin{equation}
\label{Holderbound1}
\left|\frac{1}{q_{1}}+\sum_{i=1}^{n}\frac{b_i^{1}-a_i^{1}}{q_{1}^{i+1}}-\frac{1}{q_{2}}-\sum_{i=1}^{n}\frac{b_i^{1}-a_{i}^1}{q_{2}^{i+1}}\right|
\leq C_{2} \|(x_1,y_1)-(x_2,y_2)\|_{\infty}.
\end{equation}
Here $C_2:=C_1M^2/C$.

 We now focus on removing a $|q_1-q_2|$ term from the left hand side of \eqref{Holderbound1}:
\begin{align*}
&\left|\frac{1}{q_{1}}+\sum_{i=1}^{n}\frac{b_i^{1}-a_i^{1}}{q_{1}^{i+1}}-\frac{1}{q_{2}}-\sum_{i=1}^{n}\frac{b_i^{1}-a_{i}^1}{q_{2}^{i+1}}\right|\\
=&\left|\frac{q_{2}-q_{1}}{q_{1}q_{2}}+\sum_{i=1}^{n}\frac{(b_i^{1}-a_i^{1})(q_{2}^{i+1}-q_{1}^{i+1})}{q_{1}^{i+1}q_{2}^{i+1}}\right|\\
=&\left|\frac{q_{2}-q_{1}}{q_{1}q_{2}}+(q_2-q_1)\sum_{i=1}^{n}\frac{(b_i^{1}-a_i^{1})(q_{2}^{i}+q_{2}^{i-1}q_{1}+\cdots +q_{1}^i)}{q_{1}^{i+1}q_{2}^{i+1}}\right|\\
=&|q_2-q_1|\left|\frac{1}{q_{1}q_{2}}+\sum_{i=1}^{n}\frac{(b_i^{1}-a_i^{1})(q_{2}^{i}+q_{2}^{i-1}q_{1}+\cdots +q_{1}^i)}{q_{1}^{i+1}q_{2}^{i+1}}\right|.
\end{align*}
Lemma \ref{l45qq} tells us that the second term in the product in the final line can be bounded below by a constant. Therefore if we combine Lemma \ref{l45qq} and \eqref{Holderbound1} with the above, we can assert that 
\begin{equation}
\label{Holderpart1}
|q_{2}-q_{1}|\leq C_3\|(x_1,y_1)-(x_2,y_2)\|_{\infty}.
\end{equation}
Here $C_3:=C_2/C_0$ and $C_0$ is as in Lemma \ref{l45qq}.

It remains to show that there exists $C'>0$ such that
$$|\pi_{q_1}(m(a_i^{1}))-\pi_{q_2}(m(a_i^{2}))|<C'\|(x_1,y_1)-(x_2,y_2)\|_{\infty}.$$
Note that
 $$\pi_{q_1}(m(a_i^{1}))=\frac{m}{q_1}+\sum_{i=1}^{\infty}\frac{a_i^{1}}{q_{1}^{i+1}}\quad \text{and}\quad \pi_{q_2}(m(a_i^{2}))=\frac{m}{q_2}+\sum_{i=1}^{\infty}\frac{a_i^{2}}{q_{2}^{i+1}},$$
and $a_{i}^{1}=a_{i}^{2}$ for $1\leq i\leq n$. Therefore
\begin{align}
\label{e421q}
&\left|\pi_{q_1}(m(a_i^{1}))-\pi_{q_2}(m(a_i^{2}))\right|\nonumber\\
\leq &  \left|\frac{m}{q_1}-\frac{m}{q_2}\right|+\left|\sum_{i=1}^{n}\frac{a_i^{1}}{q_{1}^{i+1}}- \sum_{i=1}^{n}\frac{a_i^{1}}{q_{2}^{i+1}}\right|+\left|\sum_{i=n+1}^{\infty}\frac{a_i^{1}}{q_{1}^{i+1}}- \sum_{i=n+1}^{\infty}\frac{a_i^{2}}{q_{2}^{i+1}}\right|\nonumber\\
\leq & \left|\frac{m}{q_1}-\frac{m}{q_2}\right|+\left|\sum_{i=1}^{n}\frac{a_{i}^1(q_{2}^{i+1}-q_{1}^{i+1})}{q_{1}^{i+1}q_{2}^{i+1}}\right|+\left|\sum_{i=n+1}^{\infty}\frac{a_i^{1}}{q_{1}^{i+1}}- \sum_{i=n+1}^{\infty}\frac{a_i^{2}}{q_{2}^{i+1}}\right|\nonumber\\
=&|q_{2}-q_{1}|\Big(\left|\frac{m}{q_1q_2}\right|+\left|\sum_{i=1}^{n}\frac{a_{i}^1(q_{2}^{i}+q_{2}^{i-1}q_{1}+\cdots+q_{1}^{i})}{q_{1}^{i+1}q_{2}^{i+1}}\right|\Big)+\left|\sum_{i=n+1}^{\infty}\frac{a_i^{1}}{q_{1}^{i+1}}- \sum_{i=n+1}^{\infty}\frac{a_i^{2}}{q_{2}^{i+1}}\right|.
\end{align}
There exists $C_{4}:= M\left(q_{KL}^{-2}+(q_{KL}-1)^{-2}\right)>0$ and $C_{5}:=M/(q_0^2-q_0)>0$ such that
\begin{equation}\label{e422q}
\begin{split}
\left|\frac{m}{q_1q_2}\right|+\left|\sum_{i=1}^{n}\frac{a_{i}^1(q_{2}^{i}+q_{2}^{i-1}q_{1}+\cdots+q_{1}^{i})}{q_{1}^{i+1}q_{2}^{i+1}}\right|&\leq M\left(\frac{1}{q_{KL}^2}+\sum_{i=1}^n\frac{i+1}{q_{KL}^{i+2}}\right)\\
&\leq M\left(\frac{1}{q_{KL}^2}+\sum_{i=0}^\infty\frac{i+1}{q_{KL}^{i+2}}\right)=C_4
\end{split}
\end{equation} and
$$\left|\sum_{i=n+1}^{\infty}\frac{a_i^{1}}{q_{1}^{i+1}}- \sum_{i=n+1}^{\infty}\frac{a_i^{2}}{q_{2}^{i+1}}\right|\leq \sum_{i=n+1}^\infty \frac{M}{q_0^{i+1}}=C_{5}q_{0}^{-n}.$$ Using the line above and \eqref{e26},
it follows that there exists $C_{6}:=C_5M^2/C>0$ such that
\begin{equation}
\label{part1}
\left|\sum_{i=n+1}^{\infty}\frac{a_i^{1}}{q_{1}^{i+1}}- \sum_{i=n+1}^{\infty}\frac{a_i^{2}}{q_{2}^{i+1}}\right|\leq C_{6}\|(x_1,y_1)-(x_2,y_2)\|_{\infty}.
\end{equation}
Using \eqref{Holderpart1}, \eqref{e422q}, and \eqref{part1}, we see that \eqref{e421q}  implies that there exists $C_{7}:=C_3C_4+C_6>0$ such that 
\begin{equation}
\label{Holderpart2}
|\pi_{q_1}(m(a_i^{1}))-\pi_{q_2}(m(a_i^{2}))|\leq C_{7}\|(x_1,y_1)-(x_2,y_2)\|_{\infty}.
\end{equation}
Equations \eqref{Holderpart1} and \eqref{Holderpart2} together imply that there exists $C':=\max\{C_3, C_7\}>0$ such that 
$$\|f(x_1,y_1)-f(x_2,y_2)\|_{\infty}\leq C'\|(x_1,y_1)-(x_2,y_2)\|_{\infty}.$$
\end{proof}

We have now proved all of the technical results we need to prove  Theorem \ref{t11}. We just require the following well known theorem due to Marstrand and a lemma which both can be found in the book by Bishop and Peres \cite{BisPer}.
	\begin{theorem}(Marstrand slicing theorem)\label{Marstrand slicing theorem}
		Let $E\subset \mathbb{R}^{2}$.
		\begin{enumerate}[\upshape (i)]
			\item If $\dim_{H}E<1$ then 
			\begin{equation}\label{e51qqqq}
			\{(q,y):y\in\mathbb{R}\}\cap E=\emptyset
			\end{equation} 
			for Lebesgue almost every $q.$ 
			\item If  $\dim_{H}E\geq 1$ then $$\dim_{H}(\{(q,y):y\in\mathbb{R}\}\cap E\})\leq \dim_{H}E-1$$ for Lebesgue almost every $q.$
		\end{enumerate}
	\end{theorem}
\begin{lemma}\label{l52}
 Suppose $K\subset \mathbb{R}^d$ and $f: \mathbb{R}^d\to \mathbb{R}^n$ is Lipschitz. Then $\dim_H f(K)\leq \dim_H K$.
\end{lemma}

\begin{proof}[Proof of Theorem \ref{t11}]
As previously remarked, to prove Theorem \ref{t11} it suffices to show that for our fixed choice of $\epsilon$, for Lebesgue almost every $q\in J$ we have
\begin{equation}\label{e51qqq}
\dim_{H}\uu_{q}^{2}\leq \max\{2\dim_H \uu_{q}-1+\epsilon, 0\}.
\end{equation}
The set $A$ satisfies $A\subseteq \bigcup_{m=0}^{M-1} (\uu_{q_0}+m)/q_0\times (\uu_{q_0}+m+1)/q_0.$ In \cite{ABBK2019} it was shown that $\dim_{H}\uu_{q_0}=\dim_{B}\uu_{q_0}$ for any $q\in (1,M+1].$ Therefore using well know properties of the Cartesian product of sets, see for instance \cite{Falconer_1990}, we know that $$\dim_{H}\left(\bigcup_{m=0}^{M-1} (\uu_{q_0}+m)/q_0\times (\uu_{q_0}+m+1)/q_0\right)=2\dim_{H}\uu_{q_0}.$$ This equality and the inclusion above imply
\begin{equation}\label{e51qq}\dim_{H}A\leq 2\dim_{H}\uu_{q_0}.\end{equation}
Since $f(A)=B$ and $f$ is Lipschitz by Proposition \ref{p24}, \eqref{e51qq} and Lemma \ref{l52} tell us that
\begin{equation}\label{e51q}
\dim_{H}B\leq 2\dim_{H}\uu_{q_0}.
\end{equation}
For any $q\in J$, it follows Lemma \ref{l43} that  \begin{equation}
	\label{e51}
\dim_{H}(\{(q,y):y\in\mathbb{R}\}\cap B)=\dim_{H}\uu_{q}^{2}.
	\end{equation}
 Theorem \ref{Marstrand slicing theorem} and \eqref{e51q} imply that for Lebesgue almost every $q\in J$ we have
 $$\dim_{H}(\{(q,y):y\in\mathbb{R}\}\cap B)\leq \max\{2\dim_{H}\uu_{q_0}-1,0\}.$$ Therefore \eqref{e51} implies that for Lebesgue almost every $q\in J$ we have
 \begin{equation}\label{e52}
 \dim_{H}\uu_{q}^{2}\leq \max\{2\dim_{H}\uu_{q_0}-1,0\}.
\end{equation} Let $q'\in J$ be an arbitrary element of the full measure subset of $J$ for which \eqref{e52} holds. Item (v) in our list of properties for the collection of intervals $\{J_i\}$ tells us that \begin{equation*}
\label{dimensionerror}
|\dim_{H}\uu_{q'}-\dim_{H}\uu_{q_0}|<\epsilon/2.
\end{equation*} Using this fact and \eqref{e52}, we see that $q'$ satisfies 
$$ \dim_{H}\uu_{q'}^{2}\leq \max\{2\dim_{H}\uu_{q'}-1+\epsilon,0\}.$$ Since $q'$ was an arbitrary element of a full measure subset of $J,$ this completes our proof that \eqref{e51qqq} holds for Lebesgue almost every $q\in J$.
\end{proof}


We now briefly explain how the argument used to prove Theorem \ref{t11} can be adapted to prove Theorem \ref{t12}.

\begin{proof}[Proof of Theorem \ref{t12}]
We recall that $$O:=\left\{q\in (q_{KL}, M+1):\dim_H\uu_q<\frac{1}{2}\right\}.$$ Just as in the proof of Theorem \ref{t12}, we can cover $O$ up to a set of measure zero by a collection of intervals $\{J_i\}$ satisfying properties (ii)-(v). We can also assume that each $J_i$ is chosen such that $q_i:=\min\{q:q\in J_i\}$ satisfies $q_i\in O$. 

We now fix an arbitrary $J$ from the collection $\{J_i\}$ and let $q_0:=\min\{q:q\in J\}$. To prove Theorem \ref{t12}, it suffices to show that for Lebesgue almost every $q\in J$ we have $\uu_{q}^{j}=\emptyset$ for $j\geq 2$. Moreover, because of Lemma \ref{l22qq}, it is in fact sufficient to show that for Lebesgue almost every $q\in J$ we have $\uu_{q}^{2}=\emptyset.$ We now proceed as in our proof of Theorem \ref{t11} and define the sets $A,B$ and the map $f$ in the same way. Importantly \eqref{e51q} tells us that $\dim_{H}B\leq 2\dim_{H}\uu_{q_0}$. We chose our intervals $\{J_i\}$ in the such a way that $q_i\in O$. Therefore $q_0\in O$ and we have $$\dim_{H}B<1.$$ Applying Theorem \ref{Marstrand slicing theorem} we may conclude that that $$\{(q,y):y\in \mathbb{R}\}\cap B=\emptyset$$ for Lebesgue almost every $q\in J$. Our theorem now follows via an application of Lemma \ref{l43} which tells us that if $\{(q,y):y\in\mathbb{R}\}\cap B=\emptyset$ then $\uu_{q}^2=\emptyset$.
\end{proof}

\thanks{\textbf{Acknowledgements.} The second author was supported by the National Natural Science Foundation
of China (NSFC)  \#11871348, \#61972265.}


\begin{thebibliography}{100}


\bibitem{ABBK2019} R. Alcaraz Barrera, S. Baker and D. R. Kong.
Entropy, topological transitivity, and dimensional properties of unique q-expansions.
\emph{Trans. Amer. Math. Soc.} \textbf{371}(5) (2019), 3209--3258.


\bibitem{AK2019}
P. Allaart, and D. R. Kong.
On the continuity of the Hausdorff dimension of the univoque set
\emph{Adv. Math. } \textbf{354} (2019), 106729.


\bibitem{BK2007}
C. Baiocchi and V.  Komornik.
 Greedy and quasi-greedy expansions in non-integer bases. Preprint, 2007,
arXiv: 0710.3001v1.

\bibitem{BisPer} C. Bishop, Y. Peres, \textit{Fractals in Probability and Analysis} (Cambridge Studies in Advanced Mathematics). Cambridge: Cambridge University Press. doi:10.1017/9781316460238


\bibitem{B2014}
S. Baker.
Generalized golden ratios over integer alphabets.
\emph{Integers} \textbf{14} (2014), Paper No. A15, 28 pp.

\bibitem{B2015}
S. Baker.
On small bases which admit countably many expansions.
\emph{J. Number Theory}  \textbf{147}  (2015), 515--532.

\bibitem{BS2014}
S. Baker and N. Sidorov.
Expansions in non-integer bases: lower order revisited.
\emph{Integers}  \textbf{14} (2014), Paper No. A57, 15 pp.



\bibitem{DD2007}
K. Dajani and M. de Vries.
Invariant densities for random $\beta$-expansions.
\emph{J. Eur. Math. Soc.} \textbf{9}(1) (2007), 157--176.

\bibitem{DK1995}
Z. Dar{\'o}czy  and I. K{\'a}tai.
On the structure of univoque numbers.
\emph{Publ. Math. Debrecen} \textbf{46}(3-4) (1995),  385--408.

\bibitem{DK2009}
M. de Vries and V. Komornik.
Unique expansions of real numbers,
\emph{Adv. Math.}  \textbf{221}(2) (2009), 390--427.


\bibitem{DKL2016}
M. de Vries,  V. Komornik and P. Loreti.
Topology of the set of univoque bases.
\emph{Topology Appl.}  \textbf{205} (2016), 117--137.

\bibitem{EHJ1991}
P. Erd\H os, M. Horv\'ath and I. Jo\'o.
On the uniqueness of the expansions $1=\sum q^{-n_i}$.
\emph{Acta Math. Hungar.} \textbf{58}(3-4) (1991), 333--342.

\bibitem{EJ1992}
P. Erd\H os  and I. Jo\'o.
On the number of expansions {$1=\sum  q^{-n_i}$}.
\emph{Ann. Univ. Sci. Budapest. E\"otv\"os Sect. Math. } \textbf{35} (1992), 129--132.

\bibitem{EJK1990}
P. Erd\H os, I. Jo\'o and V. Komornik.
 Characterization of the unique expansions $1=\sum_{i=1}^\infty q^{-n_i}$ and related problems.
\emph{Bull. Soc. Math. France}  \textbf{118}(3) (1990), 377--390.

\bibitem{Falconer_1990}
K. Falconer, K.
\emph{Fractal Geometry.
Mathematical Foundations and Applications.}
John Wiley \& Sons Ltd., Chichester, 1990.

\bibitem{GS2001}
P. Glendinning and N. Sidorov.
Unique representations of real numbers in non-integer bases.
\emph{Math. Res. Lett.} \textbf{8}(4) (2001), 535--543.



\bibitem{Kallos1999}
G. Kall\'{o}s.
The structure of the univoque set in the small case.
\emph{Publ. Math. Debrecen } \textbf{54}(1-2) (1999), 153--164.

\bibitem{Kallos2001}
G. Kall\'{o}s.
The structure of the univoque set in the big case.
\emph{Publ. Math. Debrecen}  \textbf{59}(3-4) (2001), 471--489.




\bibitem{KK2018}
V. Komornik and D. R. Kong.
Bases in which some numbers have exactly two expansions.
\emph{J. Number Theory} \textbf{195} (2019), 226--268.

\bibitem{KKL2017}
V. Komornik, D. R. Kong and W. X. Li.
Hausdorff dimension of univoque sets and devil's staircase.
\emph{Adv. Math. } \textbf{305} (2017), 165--196.


\bibitem{KL1998}
V. Komornik and P. Loreti.
Unique developments in non-integer bases.
\emph{Amer. Math. Monthly}  \textbf{105}(7) (1998), 636--639.

\bibitem{KL2002}
V. Komornik and P. Loreti.
Subexpansions, superexpansions and uniqueness properties in noninteger bases.
\emph{Period. Math. Hungar. } \textbf{44}(2) (2002) , 197--218.


\bibitem{KL2015}
D. R. Kong and W. X.  Li.
Hausdorff dimension of unique beta  expansions.
\emph{Nonlinearity }  \textbf{28}(1) (2015), 187--209.

\bibitem{KLD2010}
D. R. Kong, W. X.  Li and F. M. Dekking.
Intersections of homogeneous Cantor sets and beta-expansions.
\emph{Nonlinearity}  \textbf{23}(11) (2010), 2815--2834.

\bibitem{KLZ2017}
D. R. Kong, W. X. Li and Y. R. Zou.
On small bases which admit points with two expansions.
J. Number Theory  \textbf{173} (2017), 100--128.


\bibitem{Parry 1960}
W. Parry. \emph{On the $\beta$-expansions of real numbers.} Acta Math. Acad. Sci. Hungar. \textbf{11} (1960), 401--416.
\bibitem{PerSol} Y. Peres, B. Solomyak, \textit{Absolute continuity of Bernoulli convolutions, a simple proof,} Math. Res. Lett. \textbf{3}(2) (1996), 231--239.


\bibitem{R1957}
A. R\'enyi.
Representations for real numbers and their ergodic properties.
\emph{Acta  Math. Hungar.} \textbf{8} (1957), 477--493.

\bibitem{S2003}
N. Sidorov.
Almost every number has a continuum of  {$\beta$}-expansions.
\emph{Amer. Math. Monthly }  \textbf{110}(9) (2003), 838--842.

\bibitem{S2009}
N. Sidorov.
{Expansions in non-integer bases: lower, middle and top orders.}
\emph{J. Number Theory }  \textbf{129}(4) (2009), 741--754.

\bibitem{ZK2015}
Y. R. Zou and D. R. Kong.
On a problem of countable expansions.
\emph{J. Number Theory}  \textbf{158} (2016), 134--150.

\bibitem{ZWLB2016}
Y. R. Zou, L. J.  Wang, J. Lu and S. Baker.
On small bases for which 1 has countably many expansions.
\emph{Mathematika}  \textbf{62}(2) (2016), 362--377.

\bibitem{ZKL2021} Y. R. Zou, J. Lu, V. Komornik, Hausdorff dimension of  multiple expansions, submitted.
\end{thebibliography}
\end{document}